\documentclass{amsart}
\usepackage{amsfonts}
\input{epsf}
\usepackage{amsmath, amsthm, amssymb}
\usepackage{latexsym}
\newtheorem{theorem}{Theorem}[section]

\theoremstyle{definition}

\theoremstyle{remark}

\numberwithin{equation}{section}

\begin{document}
\title{Imbalance attractors for a strategic model of market microstructure}

\author{Ted Theodosopoulos}
\address{IKOS Research\\ 9 Castle Square, Brighton, East Sussex, BN1 1DZ, UK}
%\curraddr{}
\email{ptaetheo@earthlink.net}
%\urladdr{faculty.lebow.drexel.edu/TheodosopoulosT}
%\thanks{}
\author{Ming Yuen}
\address{Department of Mathematics 
\\ Drexel University \\ Philadelphia, PA 19104, USA}
%\curraddr{}
\email{mmy23@drexel.edu}
%\thanks{}

%    General info
%\subjclass[2000]{37B10, 82C20, 60K35}

\date{May 15, 2006.}

%\dedicatory{This paper is dedicated to our authors.}


\begin{abstract}
In this paper we extend the series of our studies on the properties of an interacting particle model for market microstructure.  In our earlier
work we defined a Markov process on the majority opinion of the agents, obtained
the transition probabilities and analyzed the martingale properties of the ensuing wealth process.  Here we relax the assumption on the choices of individual agents by allowing mixed strategies, offering opportunities for the agents to gain intermediate submartingale exposure for their individual wealth processes.  We develop a novel two-dimensional spin system to model the critical regions of the wealth process as a reflection of the agents' behaviors.  We exhibit strategic conflicts between individual market participants and the market as a whole, and identify a new source of uncertainty arising from `reinforced expectations'.
\end{abstract}

\maketitle

\section{Introduction}\label{sec:int}

A series of interacting particle models have been proposed over the past decade 
to describe the dynamics of agent-based market microstructure models \cite{aoki, bouc, horv, lux, majk}.  Often, these models use techniques from non-equilibrium statistical mechanics to investigate the tradeoff between the safety of trend following and the necessity of anticipating the majority in order to generate profits \cite{chal, cool, mars, rss}.

In our earlier work, we extended such a spin market model, originally introduced
by Bornholdt et.al. \cite{born}, and characterized analytically its invariant measure.
Also, we investigated the aggregate market wealth dynamics \cite{thmi}.  This work led
to an appreciation of the asymmetry between different quasi-stable supply-demand
imbalance levels.

The current paper extends the earlier framework by allowing the agents, with a certain
probability, to choose strategically their next action, rather than be subject to
the Hamiltonian dynamics on the space of spin configurations.  

In our revised model we extend the state space of our Markov Chain to include a set of supply/demand imbalance levels which represent the market participant's expectations.  The resulting hybrid dynamics represent a mixture of the original hamiltonian stochastic dynamics presented in \cite{theo} and the strategic dynamics of these expectation sets.  We use this setting to investigate the concept of strategic stability that we introduced in \cite{thmi} and discover situations where the individual agent's strategies are in conflict with the overall market interest.  This finding provides a setting in which to investigate the social benefits of market coordination.

Furthermore, the extended stochastic process we introduce here possesses regimes of non-unique invariant measure.  We interpret the resulting path dependence as an example of reinforced belief dynamics, whereby an arbitrary choice attains physical importance purely because a sufficient number of agents share it.  Such `self fulfilling prophecies' are a hallmark of the strongly interacting nature of markets, and constitute a main ingredient of any behavioral market model.  It is our hope that the current model can help bridge the conceptual gap between physically motivated spin market models and their behavioral counterparts.

\section{A Two-Dimensional Spin Process}\label{sec:2d}

In this study, we modify the spin market model presented in \cite{theo}.  Specifically, we will allow each agent to act strategically (in a sense to be defined below) with a fixed exogenous probability $q$.  Let $X$ denote the set of spin configurations on a lattice on the $d$-dimensional torus\footnote{We use the notation ${\mathcal T}^d$ to denote the object $\underbrace{ {\mathcal S}^1 \times \ldots \times {\mathcal S}^1}_d$.} $Y \doteq \left({\mathcal Z}/L \right)^d \subset {\mathcal T}^d$, i.e. $X \subset \{-1,1 \}^Y$, for an appropriately chosen $L$ so that $|Y|=N$.  It is convenient to identify the elements of $Y$ with unique labels from the set $\{1, 2, \ldots, N\}$.  Formally, we accomplish this by choosing a homeomorphism $\varphi: Y \longrightarrow \{1, 2, \ldots, N\}$.
%and identify $Y$ with its image under $\varphi$.  By this abuse of notation, $y %\in Y$ will stand for the agent at location $y$ as well as the label %$\varphi(y)$.  

The state space of our revised process is ${\mathcal X} = X \times X$ where the second component of a typical state $(\zeta_1, \zeta_2)$ represents a subset $G \subset Y$, i.e. $y \in G \Longleftrightarrow \zeta_2(y) = 1$.  This can be seen as the assignment of a `two-dimensional spin' to the lattice $Y$, i.e. ${\mathcal X} \ni \vec{\zeta} : Y \longrightarrow \{-1, 1\}^2$.

The path of a typical element of ${\mathcal X}$ is given by $\eta: Y \times ( 0, \infty ) \longrightarrow \{-1,1\}^2$ and each site $x \in Y$ is endowed with a (typically $\ell_1$) neighborhood ${\mathcal N} (x) \subset Y$ it inherits from the natural topology on the torus ${\mathcal T}^d$.  As before, for each $x \in Y$, ${\mathcal N} (x)$ is a uniformly chosen random subset of $Y$, of cardinality $2d$.  To be more precise, for a set $A$ and a positive integer $k$, let $F(A,k) = \left\{ \left(a_1, a_2, \ldots, a_k \right) \in A^k | a_i \neq a_j \mbox{         } \forall i,j=1,2,\ldots,k \right\}$.  Then, for any $x \in Y$, let $\left\{ {\mathcal N} (x, \cdot) \right\}$ be a family of iid uniform random variables taking values in $F \left( Y \setminus \{x\}, 2d \right)$.

We construct a continuous time Markov process with transitions occurring at exponentially distributed epochs, $T_n$, with rate 1.  We proceed as before to construct a transition matrix for the spins, based on the following interaction potential:
$$h(x,T_n) = \sum_{y \in {\mathcal N}(x,n)} \eta_1 (y,T_n) - \alpha \eta_1 (x,T_n) N^{-1} \left|\sum_{y \in Y} \eta_1 (y,T_n) \right|,$$
where $\alpha>0$ is the coupling constant between local and global interactions.  As discussed in \cite{theo}, observe that, as a result of the neighborhood randomization procedure we introduced, $h(x,t)$ depends only on the first dimension of the spin at $x$, $\eta_1(x,t)$, and $N^+(t) \doteq \sum_{x \in Y} \left[\eta_1 (x, t) \vee 0 \right]$, i.e. the total number of positive spins.  Specifically, consider two independent families of iid hypergeometric random variables, $\{U_1, U_2, \ldots \}$ and $\{V_1, V_2, \ldots \}$, with parameter vectors $(N-1, i-1, 2d)$ and $(N-1, i, 2d)$ respectively.  Then, for $t \in \left. \left[T_n, T_{n+1} \right. \right)$, define
\begin{eqnarray*}
h_{+1}(i, t) & = & 2(U_n -d) - \alpha \left| {\frac {2i}{N}} -1 \right| \\
h_{-1}(i, t) & = & 2(V_n -d) + \alpha \left| {\frac {2i}{N}} -1 \right|.
\end{eqnarray*}
We see that, conditional on $\{N^+(T_{n-1}) = i\} \cap \{\eta_1 (x, T_{n-1}) = \pm 1 \}$, the distribution of $h(x,T_n)$ is the same as that for $h_{\pm 1}(i, T_n)$.

The flips of the first spin dimension are governed by: 
$$\rho_{ab} (q,\beta, i, t) = \left[{\frac {ai}{N}} +{\frac {1-a}{2}} \right] \left\{ {\frac {q}{1+ \exp[-2b \beta h_a (i, t)]}} + {\frac {(1-q) \left[1+ b \eta_2 (\varphi^{-1} (i), t) \right]}{2}} \right\},$$
where $(a,b) \in \{-1,1\}^2$ and $\beta$ is the normalized inverse temperature.  
At time $T_n$ (i.e. the $n^{\rm th}$ epoch) a random site $J_n$ is chosen uniformly from $Y$.  The variables $\rho_{ab}$ are interpreted as the resulting flip probabilities for the first spin dimension of $J_n$, i.e.
$$\rho_{ab} = {\bf Pr} \left(\eta_1 (J_n, T_{n-1}) = a \mbox{      } \& \mbox{      } \eta_1 (J_n, T_n) = b \left| {\mathcal B}_{n-1} \cap \left\{N^+ (T_{n-1}) = i \right\} \right. \right).$$
If we let 
$$\bar{X}_n  = {\frac {\eta_1 (J_n, T_n) - \eta_1 (J_n, T_{n-1})}{2}} = N^+ (T_n) -N^+ (T_{n-1}) \in \{-1,0,1\},$$
we have
$${\bf Pr} \left(\bar{X}_n = a \left| {\mathcal B}_{n-1} \cap \left\{N^+ (T_{n-1}) = i \right\} \right. \right) = \rho_{-a, a} (q,\beta, i, T_n).$$

As stated so far, the transition probabilities of $\bar{X}$ (and therefore $N^+$) could be random.  In particular, for $t \in \left[ \left. T_{n-1}, T_n \right) \right.$, even though $h_a (i,t)$ is measurable with respect to ${\mathcal B}_{n-1}$, $\eta_2 (\varphi^{-1} (i), t)$ isn't generically\footnote{This is because, for $t \in \left[ \left. T_{n-1}, T_n \right) \right.$, the distribution of $h_a (i,t)$ is determined by $N^+(T_{n-1})$, while, for $\beta < \infty$, the distribution of $\rho_{ab}$ (and therefore $P_{ab}$ and ${\mathcal E}_a$) depend also on whether an agent with Hamiltonian $h_a$ has spin $+1$ or $-1$, which is a Bernoulli random variable}.  In order to describe the transition probabilities for the second component of our stochastic process, we need to define the wealth process \cite{thmi}:
$$W(y,T_n) \doteq \left\{
\begin{array}{ll}
W(y,T_{n-1}) & \mbox{for $J_{n}=y$;} \\
W(y,T_{n-1}) + \eta_1 (y,T_{n-1})\Delta P(T_n) & \mbox{otherwise.}\\
\end{array}
\right.,$$
where $W(y,0)= K(y)$, representing the initial capital available to agent $y$, and $\Delta P(T_n) \doteq P(T_n) -P(T_{n-1})$ is the price change at the $n^{\rm th}$ epoch.  In what follows, we will allow a more general price impact function than in our earlier models, as follows:
$$\Delta P(T_n) = P (T_{n-1}) f \left( \bar{X}_n, N \right),$$
where $f$ is a function that maps transient supply/demand imbalances $\bar{X}_n$ to price shocks.  The only assumptions we make about this impact function $f$ is that it is nondecreasing in the first argument and nonincreasing in the second argument, and that it maps the balanced market to the fundamental price, i.e. $f(0,N) = 0$.  In what follows, we use $\gamma$ to denote the ratio of the price impact of an incremental seller to that of an incremental buyer, i.e. $\gamma = {\frac {f(-1,N)}{f(1,N)}} \leq 0$, suppressing its dependence on $N$ for convenience.

Let $P^\beta_{ab} (i)$ be independent random variables, identically distributed with the iid family $\rho_{ab} (1, \beta, i, \cdot)$.  These represent the purely stochastic transitions of the first spin dimension investigated in \cite{theo}.  Using them, we define
$${\mathcal E}^\beta_a (i) = P^\beta_{-a,a} (i) +\gamma P^\beta_{a,-a} (i),$$
which is a random variable so long as $\beta< \infty$, inheriting the randomness of $h_{\pm 1} (i, \cdot)$.  On the other hand, as discussed in \cite{theo}, in the frozen phase of this process, almost surely,
\begin{equation}
{\mathcal E}_+ (i) \doteq \lim_{\beta \rightarrow \infty} {\mathcal E}^\beta_+ (i) = \left( 1 -{\frac {i}{N}} \right) {\bf Pr} \left( h_{-1} (i, \cdot)>0 \right) + {\frac {\gamma i}{N}} {\bf Pr} \left( h_{+1} (i, \cdot)<0 \right), \label{eq:frozen}
\end{equation}
which is not random.

The role that $\rho$ played for the transitions of the first spin dimension will now be played by $\lambda$ for the second spin dimension.  Specifically, for $(a,b) \in \{-1,1\}^2$ and $t \in \left[ \left. T_{n-1}, T_n \right) \right.$, let
$$\lambda_{ab}(q, \beta, x, t) = {\bf Pr} \left( \left. a{\mathcal E}^\beta_+ (\varphi(x)) +ab \left(1 -{\frac {1}{q}} \right) \left(1 -{\frac {\varphi(x)}{N}} \right) > 0 \right| {\mathcal B}_{n-1} \cap \left\{ \eta_2 (x, T_{n-1}) = a \right\} \right).$$
Note that the event $ {\mathcal E}^\beta_+ (i) = -b \left(1 -{\frac {1}{q}} \right) \left(1 -{\frac {i}{N}} \right)$, which occurs generically with probability zero, would be added to the computation of $\lambda_{aa}$, i.e. the probability of not flipping the second spin dimension.  The variables $\lambda_{ab}$ are interpreted as the flip probabilities for the second spin dimension of at the $n^{\rm th}$ epoch, i.e. for $t \in \left[ \left. T_{n-1}, T_n \right) \right.$ and all $x \in Y$,
$$\lambda_{ab} (q, \beta, x, t) = {\bf Pr} \left(\eta_2 (x, T_n) = b \left| {\mathcal B}_{n-1} \cap \left\{\eta_2 (x, T_{n-1}) = a \right\} \right. \right).$$

\section{Convergence Results}\label{sec:convergence}

Observe that, unlike the first spin dimension, the second one is updated synchronously for all sites.  Furthermore, using (\ref{eq:frozen}) we see that in the frozen phase, $\lambda_{ab} \in \{0,1\}$, because ${\mathcal E}_+^\beta$ is no longer random.  In particular, let 
$${\mathcal B}= \left\{x \in Y |{\mathcal E}_+^\beta \left(\varphi (x) \right) \geq \left( 1- {\frac {1}{q}} \right) \left(1- {\frac {\varphi(x)}{N}} \right) \right\}$$ and 
$${\mathcal C}= \left\{x \in Y |{\mathcal E}_+^\beta \left(\varphi (x) \right) \leq \left( 1- {\frac {1}{q}} \right) {\frac {\varphi(x)}{N}} \right\},$$ 
$${\mathcal A}_1 = {\mathcal B} \cap {\mathcal C}^c \mbox{ , } {\mathcal A}_2 = {\mathcal B} \cap {\mathcal C} \mbox{ , } {\mathcal A}_3 = {\mathcal B}^c \cap {\mathcal C}^c \mbox{ and } {\mathcal A}_4 = {\mathcal B}^c \cap {\mathcal C}.$$  
Then:

\begin{theorem}
\label{theo1}
There are four distinct behaviors for $\eta_2$:
\begin{eqnarray*}
& & x \in \varphi^{-1} \left( {\mathcal A}_1 \right) \Longrightarrow \eta_2 (x,t)=1 \mbox{ for } t \geq T_1 \\
& & x \in \varphi^{-1} \left( {\mathcal A}_4 \right) \Longrightarrow \eta_2 (x,t)=-1 \mbox { for } t \geq T_1 \\
& & x \in \varphi^{-1} \left( {\mathcal A}_2 \right) \Longrightarrow \eta_2 (x,t)=\eta_2 (x,T_0) \mbox{ for } t \geq T_0 \\
& & x \in \varphi^{-1} \left( {\mathcal A}_3 \right) \Longrightarrow \eta_2 (x,t)=-\eta_2 (x,T_{n-1}) \mbox{ for } t \in \left[ \left. T_n, T_{n+1} \right) \right. \mbox{ and } n>0.
\end{eqnarray*}
The third case entails path dependence, while the fourth leads to oscillating behavior.
\end{theorem}

In the rest of this paper, we will focus on the frozen phase ($\beta \rightarrow \infty$).  Let $P_{ab} (i) = \lim_{\beta \rightarrow \infty}  P_{ab}^\beta (i)$.  As we saw in \cite{theo},
\begin{eqnarray*}
& & \begin{array}{cc} \bar{P}_{\scriptscriptstyle ++}(i) = \bar{P}_{\scriptscriptstyle --}(i) = 0 & \;\;\;\;\; \mbox{if $i \not\in \left[N \left( {\frac {1}{2}} - {\frac {d}{\alpha}} \right), N \left( {\frac {1}{2}} + {\frac {d}{\alpha}} \right) \right]$} \end{array} \\
& & \begin{array}{cc} \bar{P}_{\scriptscriptstyle ++}(i) = {\displaystyle \sum_{\scriptscriptstyle j=(d + c) \vee (i+2d-N)}^{\scriptscriptstyle 2d \wedge i}} f_+(i,j) & \;\;\;\;\; \mbox{if $i \in \left[N \left( {\frac {1}{2}} - {\frac {d}{\alpha}} \right), N \left( {\frac {1}{2}} + {\frac {d}{\alpha}} \right) \right]$} \end{array} \\
& & \begin{array}{cc} \bar{P}_{\scriptscriptstyle --}(i) = {\displaystyle \sum_{\scriptscriptstyle j=0 \vee (i+2d-N)}^{\scriptscriptstyle (d-c) \wedge i}} f_-(i,j) & \;\;\;\;\; \mbox{ if $i \in \left[N \left( {\frac {1}{2}} - {\frac {d}{\alpha}} \right), N \left( {\frac {1}{2}} + {\frac {d}{\alpha}} \right) \right]$} \end{array} \\
& & P_{\scriptscriptstyle ++} = {\frac {i}{N}} \bar{P}_{\scriptscriptstyle ++}, \;\;\;\;\; P_{\scriptscriptstyle +-} = {\frac {i}{N}} \left( 1 -\bar{P}_{\scriptscriptstyle ++} \right), \;\;\;\;\; P_{\scriptscriptstyle --} = \left(1 -{\frac {i}{N}} \right) \bar{P}_{\scriptscriptstyle --}, \\
& & P_{\scriptscriptstyle -+} = \left(1 -{\frac {i}{N}} \right) \left(1 -\bar{P}_{\scriptscriptstyle --} \right)
\end{eqnarray*}
with
\begin{eqnarray*}
& & c = \left\lceil \alpha \left| {\frac {i}{N}} - {\frac {1}{2}} \right| \right\rceil , \;\;\;\;\; f_+(i,j) = C_{\scriptscriptstyle j}^{\scriptscriptstyle i-1} C_{\scriptscriptstyle 2d-j}^{\scriptscriptstyle N-i} \left/ C_{\scriptscriptstyle 2d}^{\scriptscriptstyle N-1} \right.  \;\;\;\;\; f_-(i,j) = C_{\scriptscriptstyle j}^{\scriptscriptstyle i} C_{\scriptscriptstyle 2d-j}^{\scriptscriptstyle N-i-1} \left/ C_{\scriptscriptstyle 2d}^{\scriptscriptstyle N-1} \right. \\
\end{eqnarray*}
and $C_k^n$ denoting the combinations $n$ choose $k$.

Let 
$$\pi_\infty (\ell) = \lim_{\beta \rightarrow \infty} \lim_{n \rightarrow \infty} {\rm Pr} \left( N^+ \left( T_n \right) = \ell \right).$$  
The following theorem describes the invariant measure of $N^+$ in this new two-dimensional spin process, extending Proposition 1 from \cite{theo}:

\begin{theorem}
\label{theo2}
When ${\mathcal A}_2 \cup {\mathcal A}_3 = \emptyset$ there is a unique invariant measure for $N^+$ given by:
\begin{eqnarray*}
g(\ell) & = & C^N_\ell \prod^{\ell-1}_{j=0} {\frac {1 -q\bar{P}_{\scriptscriptstyle --}(j) -(1-q) \chi_{{\mathcal A}_4}(j)}{1 -q\bar{P}_{\scriptscriptstyle ++}(j+1) -(1-q) \chi_{{\mathcal A}_1}(j+1)}} \;\;\;\;\; \mbox{   if $\ell>0$ and $g(0)=1$} \\
Z(N) & = & \sum^N_{i=1} g(i) \\
\pi_\infty (\ell) & = & \left(1+ Z(N) \right)^{-1} g(\ell) \;\;\;\;\; \mbox{   if $0 \leq \ell \leq N$ and $0$ otherwise} 
\end{eqnarray*}

On the other hand, when ${\mathcal A}_3 = \emptyset$ but ${\mathcal A}_2 \ne \emptyset$, the invariant measure described above is not unique, because $\eta_2 ( \cdot, \infty)$ on ${\mathcal A}_2$ depends on $\eta_2 (\cdot, 0)$.  In this case, any measure of the form above, substituting ${\mathcal A}_4$ with ${\mathcal A}_4 \cup \{\eta_2(\cdot,0)=-1\}$ and ${\mathcal A}_1$ with ${\mathcal A}_1 \cup \{\eta_2(\cdot,0)=+1\}$, will be invariant under the two-dimensional spin process.  Thus there are $2^{|{\mathcal A}_2|}$ different invariant measures, one for each assignment of $\pm 1$ to the second spin dimension of the points in ${\mathcal A}_2$.

Finally, when ${\mathcal A}_3 \ne \emptyset$, no invariant measure exists for the two-dimensional spin process, because $\eta_2$ never settles for the points in ${\mathcal A}_3$.
\end{theorem}

\begin{figure}
\epsfxsize=5in
\epsfbox{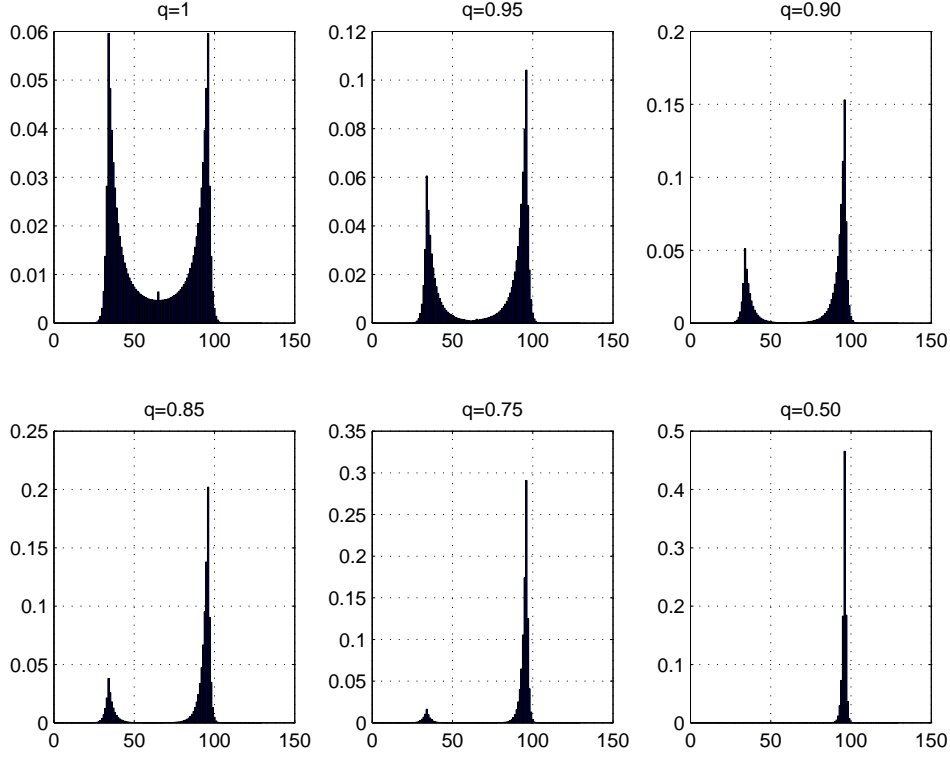}
\caption{Skewing of the invariant measure as a result of the inclusion of $q \ne 1$ ($\gamma = -0.7$, $\alpha = 4.1$, $N=128$ and $d=2$).}
\label{fig:invariant}
\end{figure}

Figure \ref{fig:invariant} shows how the invariant measure varies as we decrease $q$ away from $1$.  The case $q=1$, which was treated in \cite{theo}, always led to symmetric tri-modal distributions.  Here we observe that even a modest amount of strategic interaction leads to pronounced skew in the invariant measure.

\begin{figure}
\epsfxsize=4in
\epsfbox{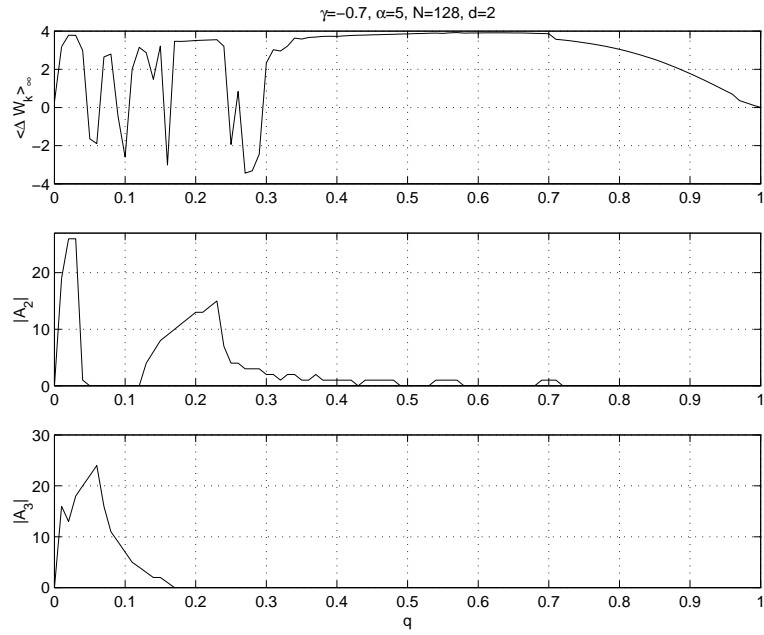}
\caption{Existence, uniqueness and martingale property of the invariant measure as a function of $q$.}
\label{fig:properties}
\end{figure}

The bottom two frames of figure \ref{fig:properties} show how the existence and uniqueness of the invariant measure described in Theorem \ref{theo2} above depend on the value of $q$ (we will return to the top frame of figure \ref{fig:properties} in section \ref{sec:wealth}).  Specifically we see that there is typically a contiguous range of small enough values of $q$ for which the invariant measure does not exist.  Thus, when strategic interactions dominate the market, agents' expectations of desirable levels of supply/demand imbalance never settle down, indefinitely forcing the agents to chase unattainable majority/minority levels.  

On the other hand, there are islands of non-uniqueness throughout the range of $q$ values, brought on by intermittent occupation of ${\mathcal A}_2$.  Except for sufficiently low values of $q$, the non-uniqueness is driven by very few (often only one) imbalance level whose status as desirable (or not) is determined by the initial conditions.  The resulting path dependence propagates throughout the invariant measure.  

\begin{figure}
\epsfxsize=4in
\epsfbox{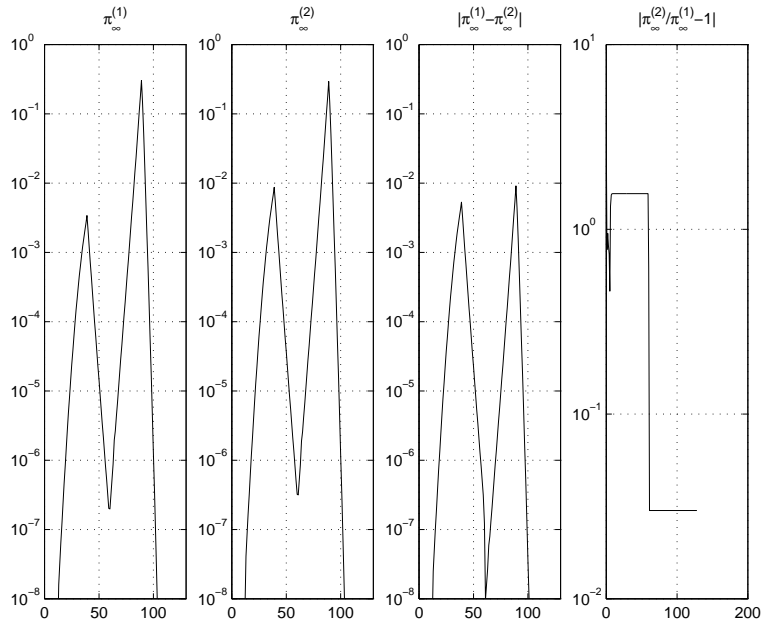}
\caption{The two invariant measures for $\gamma=-0.7$, $\alpha=5$, $N=128$, $d=2$ and $q=0.7$.  The non-uniqueness is driven by $\eta_2 \left(\varphi^{-1} (61,0) \right) = \pm 1$.}
\label{fig:nonunique}
\end{figure}

For instance, as we can see in figure \ref{fig:properties}, the case $q=0.7$ leads to a non-unique invariant measure because there is one imbalance level (it turns out to be $N^+ = 61$) which is mapped by $\eta_2$ in ${\mathcal A}_2$.  The two resulting invariant measures are shown in figure \ref{fig:nonunique}.  As we see, the flip of a single $\eta_2$ value in the initial conditions leads to more than double the limiting probability for a substantial range of seller-dominated imbalance levels.  In fact, all imbalance levels are impacted by at least a few percentage points relative change in the limiting distribution.  This phenomenon can be interpreted as a sort of fad dynamic, a positive feedback of reinforced expectations with a surprisingly global and persistent effect.

\begin{figure}
\epsfxsize=4in
\epsfbox{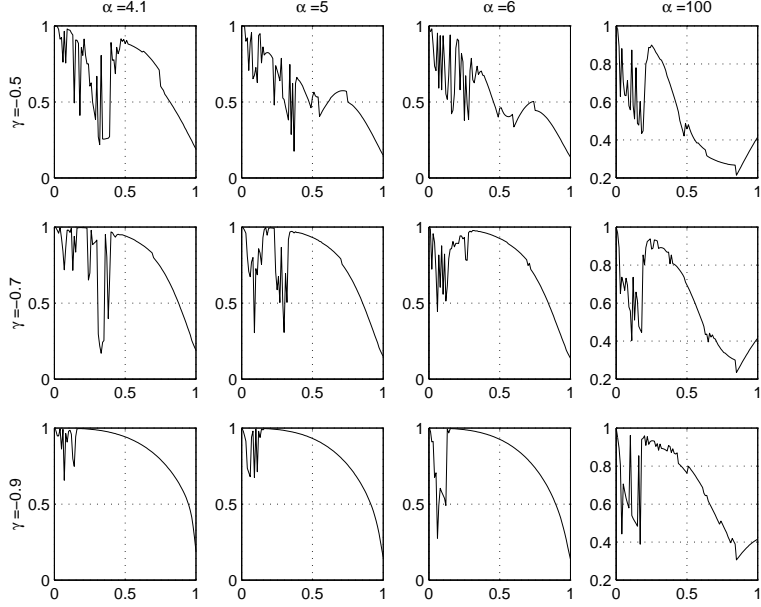}
\caption{Dependence of $\pi_\infty \left( [i^\ast -2, i^\ast +2] \right)$ on $\alpha$ and $\gamma$, where $i^\ast = \arg \max_{i \in [0,N]} \pi_\infty(i)$.}
\label{fig:rangetaila}
\end{figure}

The invariant measure changes qualitatively as $\gamma$ decreases towards $-1$ and $\alpha$ increases.  Specifically, more symmetric price impact ($\gamma \searrow -1$) leads to a steeper concentration of measure as $q$ decreases away from $1$.  In particular, as we increase the strategic interaction of the agents, all but one of the modes decrease rapidly, effectively concentrating the limiting distribution around the global mode.  This is shown in figure \ref{fig:rangetaila}, which depicts the concentration of limiting probability within two imbalance levels of the global mode.  

\begin{figure}
\epsfxsize=4in
\epsfbox{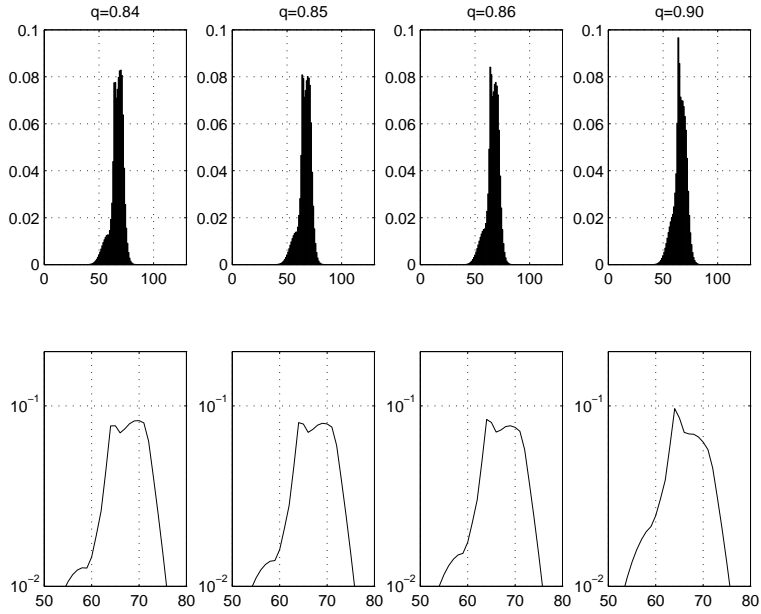}
\caption{Disappearance of bimodality as $q \nearrow 1$ for large $\alpha$ ($\gamma =-0.9$, $\alpha =100$, $N=128$, $d=2$).  The bottom frame is a blowup around the global mode of the invariant measure.}
\label{fig:unimodal}
\end{figure}

On the other hand, increasing the overall coupling constant ($\alpha \nearrow \infty$) moves the modes closer together without attenuating their height, until they are eventually merged.  The transition to unimodal limiting distribution occurs at the `kink' in the right most column of graphs in figure \ref{fig:rangetaila}, for the value of $q$ which minimizes the concentration of measure phenomenon.  This effect is shown in figure \ref{fig:unimodal}.  

Multimodal invariant measures imply the existence of distinct imbalance attractors.  Invariant measures more concentrated around their global mode imply markets with narrower price fluctuations.  Therefore, even modest increases in the amount of strategic interaction among the agents quickly decrease the spread and complexity of the resulting price fluctuations.  Moreover, this effect becomes more pronounced for more symmetric impact functions and more loosely coupled (less frustrated) agents.

\begin{figure}
\epsfxsize=3in
\epsfbox{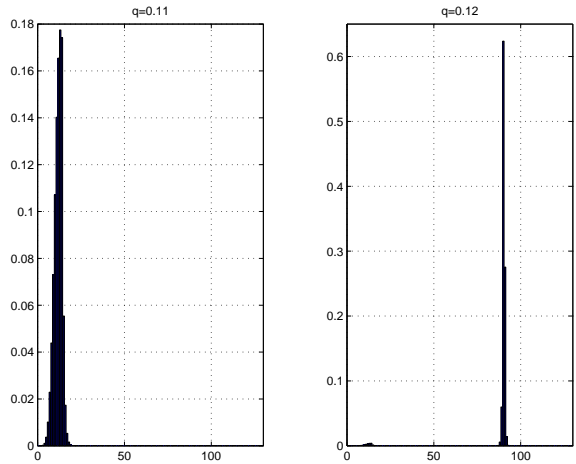}
\caption{Discontinuous transition of invariant measure for modest change in the proportion of strategic interaction $q$ ($\alpha=5, \gamma=-0.9, N=128, d=2$).}
\label{fig:jump}
\end{figure}

There also exist situations in which the invariant measure exhibits exceptional sensitivity to the amount of strategic interaction among the agents.  Figure \ref{fig:jump} shows an example of this phenomenon.  This radical shift in the imbalance attractor, from an average value of $10.84$ at $q=0.11$ to an average value of $87.71$ at $q=0.12$!  Such large scale rebalancing of the invariant measure is likely to lead to discontinuous jumps in the resulting price process, due solely to a slight decrease in the strategic component of the interaction potential.

\section{The Wealth Process Revisited}\label{sec:wealth}

We are now interested in characterizing the martingale properties of $W(y, \cdot)$ as well as the aggregate wealth, $W$, defined as  
\begin{eqnarray*}
W(T_k)&=& W(T_{k-1})+\Delta P(T_k) \sum_{y \in Y} \eta_1 (y,T_{k-1}) \\
& = & W(T_{k-1})+\Delta P(T_k) \left[2N^+(T_{k-1})-N+\bar{X}_k \right]. 
\end{eqnarray*}
Using (\ref{eq:frozen}), after some algebra we arrive at
\begin{eqnarray}
& & {\bf E} \left[W (y, T_k) -W(y, T_{k-1}) \left| {\mathcal B}_{k-1} \right. \right] = \label{eq:agentwealth} \\
& = & -\eta_1 (y, T_{k-1}) f(1,N) P (T_{k-1}) \left(1- N^{-1} \right) \left( 1+\gamma^{-1} \right){\mathcal E}_+ \left( N^+ T_{k-1} \right) \nonumber
\end{eqnarray}
and
\begin{eqnarray}
& & {\bf E} \left[W (T_k) -W(T_{k-1}) \left| {\mathcal B}_{k-1} \right. \right] = -f(1,N) P (T_{k-1}) \left( 1+\gamma^{-1} \right) \label{eq:marketwealth} \\
& & \left\{ \left[2N^+ (T_{k-1}) -N+1 \right] P_{\scriptscriptstyle -+} + \gamma \left[2N^+ (T_{k-1}) -N-1 \right] P_{\scriptscriptstyle +-} \right\}. \nonumber
\end{eqnarray}
Using the invariant measure from Theorem \ref{theo2} we can define $\langle \Delta W_k \rangle_\infty$ as the expectation, under the invariant measure $\pi_\infty$, of the expected wealth increment at the $k^{\rm th}$ transaction time:
\begin{eqnarray}
\langle \Delta W_k \rangle_\infty & = & -f(1,N) P (T_{k-1}) \left( 1+\gamma^{-1} \right) \left(1+ Z(N) \right)^{-1} \nonumber \\
& & \sum_{\ell =0}^N g(\ell) \left\{ \left[2\ell -N+1 \right] P_{\scriptscriptstyle -+}(\ell) + \gamma \left[2\ell -N-1 \right] P_{\scriptscriptstyle +-} \right\} (\ell).  \label{eq:dw}
\end{eqnarray}
We can now return to the top frame of figure \ref{fig:properties} which depicts $\langle \Delta W_k \rangle_\infty$.  As we see, the highest value of expected market wealth increments is achieved for a level of $q$ close to $0.6$.  Specifically, the maximum value, $3.9128$, is achieved for $q=0.56$ and the invariant measure is not unique, because $\eta_2 (59) \in {\mathcal A}_2$.  On the other hand the second highest value, $3.9102$, is achieved for $q=0.62$ and the invariant measure is unique.  In fact, with the exception of narrow islands restricted to $q<0.3$, the market expects to make money under the conditions in figure \ref{fig:properties}.  

\begin{figure}
\epsfxsize=4in
\epsfbox{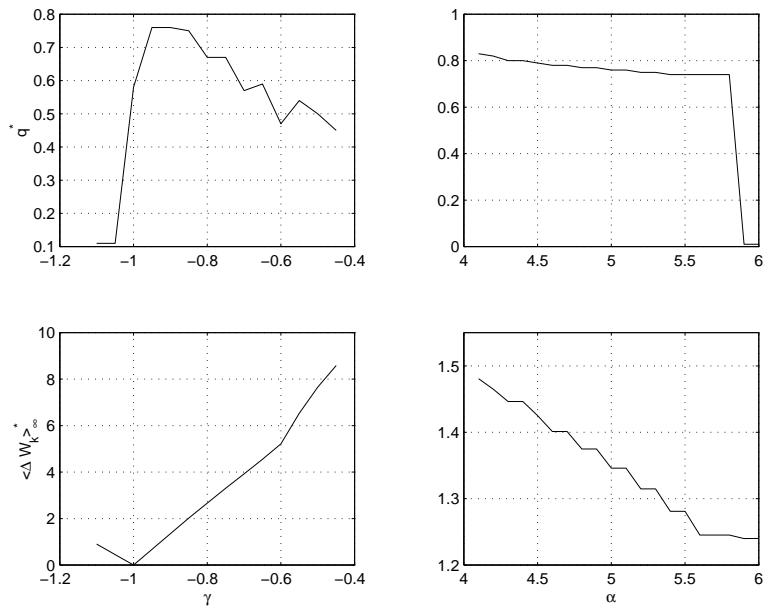}
\caption{Sensitivity of $q^\ast$ and resulting $\langle \Delta W_k \rangle_\infty^\ast$ to changes in $\alpha$ (left) and $\gamma$ (right) ($\alpha=5$ for the left panel, $\gamma=-0.9$ for the right panel, $N=128$ and $d=2$ throughout).}
\label{fig:qstrsensitivity}
\end{figure}

We proceed to iterate this computation of the optimal level of strategic interaction $q^\ast$ and the resulting maximum attainable expected market wealth increment $\langle \Delta W_k \rangle_\infty^\ast$ as we vary $\alpha$ and $\gamma$ (figure \ref{fig:qstrsensitivity}).  We can make the following observations regarding the sensitivity of the expected market wealth increments to changes in the impact function asymmetry ($\gamma$) and the overall coupling constant ($\alpha$):

\begin{itemize}
\item For all values of $\gamma$, the optimal choice of Strategic interaction level $q$ leads to non-negative expected wealth increments for the market at large.  
\item When the impact function is symmetric ($\gamma =-1$), the market can no longer make positive returns.  Thus it is the impact function asymmetry that sustains positive expectations of broad market returns.
\item Modest positive impact asymmetry (incremental buy has slightly more effect than incremental sell) leads to the minimum desirable amount of strategic interaction.
\item As the global coupling constant is increased, more strategic interaction is desirable while the attainable returns decrease.
\item When the coupling constant is increased beyond a sufficiently high value (around $5.9$ for the example shown in figure \ref{fig:qstrsensitivity}), all agents act exclusively strategically.
\end{itemize}

Consider now the majority opinion of the agents,
$$M \left(T_k \right) = \sum_{y \in Y} {\rm sgn} \left({\bf E} \left[W (y, T_k) -W(y, T_{k-1}) \left| {\mathcal B}_{k-1} \right. \right]\right).$$
Combining (\ref{eq:agentwealth}) and (\ref{eq:marketwealth}) we arrive at the following observation:

\begin{theorem}
\label{theo3}
When the market experiences a supermartingale, so does the majority of the agents.  There exist situations when the market disagrees with the majority of the agents, in the sense that their respective expected wealth increments have opposite signs.  When the market disagrees with the majority, the market is experiencing a submartingale.
\end{theorem}

Thus, there are cases when coordination of the agents to achieve a `common objective' is desirable over individual optimization of their respective strategic interaction strengths.

\section{Conclusions and Next Steps}

In this paper we presented a two-dimensional spin market model.  The first spin dimension represents the agents' decision to buy or sell, as in earlier models.  The second spin dimension encodes a subset of supply/demand imbalance levels which are deemed to be desirable, because they lead the agent under consideration to experience a submartingale wealth process.

We extended the analysis from \cite{theo} to describe the invariant measure of this novel process.  Our principal findings are that, under certain conditions the resulting invariant measure is not unique, depending sensitively on initial conditions.  We interpret this behavior as a type of `fad dynamic'.  Also, there exist conditions under which the market disagrees with the majority of the agents in its assessment of the desirability of certain imbalance levels.  Thus, coordination across the agents adds value beyond that attainable by individual optimization by the agents.

It is natural to extend further our model by allowing each agent to adjust their degree of strategic interaction ($q$) dynamically as they experience the process.  The resulting evolutionary game could provide insight into periodic imbalance attractors.  On the other hand, the thermodynamic limit ($N \rightarrow \infty$) of this process is likely to provide estimates of the occupation density around individual imbalance levels.

\bibliographystyle{amsalpha}

\end{document}